\documentclass[letterpaper,conference]{IEEEtran}
\IEEEoverridecommandlockouts
\usepackage{graphicx,epsfig,psfrag,rotate,xcolor}

\let\labelindent\relax

\usepackage{amsmath,amsfonts,amssymb,latexsym,amsthm}
\usepackage{setspace,enumerate,ifthen,subfig}
\usepackage{hyperref}
\usepackage{cite}
\usepackage{ifpdf}
\usepackage{subfiles}
\usepackage{enumitem}
\usepackage{tabularx}
\usepackage{supertabular}

\usepackage{algorithm}
\usepackage[noend]{algorithmic}
\usepackage{csquotes}
\usepackage{soul} 

\usepackage{tikz}
\usetikzlibrary{shapes,arrows}
\usetikzlibrary{shapes.geometric}
\tikzstyle{line} = [draw, -latex']
\usepackage{tikz}
\usetikzlibrary{shapes,arrows}
\tikzstyle{decision} = [diamond, draw, fill=white, 
    text width=4em, text badly centered, node distance=3cm, inner sep=0pt]
\tikzstyle{block} = [rectangle, draw, fill=blue!20, 
    text width=8em, text centered, rounded corners, minimum height=4em]
\tikzstyle{blank} = [rectangle, fill=white, 
    text width=8em, text centered, rounded corners, minimum height=4em]
\tikzstyle{line} = [draw, -latex']
\tikzstyle{vague} = [draw, -latex']
\tikzstyle{cloud} = [draw, ellipse,fill=white, node distance=3cm,
    minimum height=2em]

\definecolor{Blu}{rgb}{.255,.41,.884} 

\newtheorem{proposition}{Proposition}[section]

\newcommand{\R}{\mathbb{R}}

\allowdisplaybreaks

\begin{document}

\title{Projections onto the Set of Feasible Inputs \\ and the Set of Feasible Solutions \thanks{Contact author: \tt{jakub.marecek@ie.ibm.com}.}}

\author{\IEEEauthorblockN{Claudio Gambella}\IEEEauthorblockA{
\textit{IBM Research -- Ireland}\\
Dublin D15, Ireland
} \and \IEEEauthorblockN{Jakub Marecek} \IEEEauthorblockA{
\textit{IBM Research -- Ireland}\\
Dublin D15, Ireland
} \and \IEEEauthorblockN{Martin Mevissen}\IEEEauthorblockA{
\textit{IBM Research -- Ireland}\\
Dublin D15, Ireland
}
}

\maketitle
\thispagestyle{empty}
\pagestyle{empty}

\begin{abstract}
We study the projection onto the set of feasible inputs and the set of feasible solutions of a polynomial optimisation problem (POP).
Our motivation is increasing the robustness of solvers for POP: Without \emph{a priori} guarantees of feasibility of a particular instance,
one should like to perform the projection onto the set of feasible inputs
prior to running a solver. 
Without a certificate of optimality,
one should like to project the output of the solver onto the set of feasible solutions subsequently.
We study the computational complexity, 
formulations, and convexifications of the projections. 
Our results are illustrated on 
IEEE test cases of Alternating Current Optimal Power Flow (ACOPF) problem. 
\end{abstract}



\section{Introduction}

Polynomial optimisation is an important branch of non-convex optimisation.
For its best-known application, the alternating-current optimal power flow (ACOPF) problem 
in the management of operations of electric power systems,
there are estimates \cite{castillo2013computational} that 
a 5 \% improvement in the solvers would amount to savings of \$26B \emph{per annum}, 
in the United States alone.
Nevertheless, even this one commercially important special case remains very challenging,
after half a century of study \cite[cf.]{carpentier1962contribution}.

There are a number of challenges:
All available solvers struggle to detect infeasibility reliably and most do not respond to infeasible instances
other than by declaring infeasibility, which requires a human intervention.
Worse, leading general-purpose non-linear programming (NLP) software packages do not guarantee to provide a feasible solution, 
even in case the instance is feasible,
and often fail to provide a feasible solution in practice.
Specialised solvers for the ACOPF based on ideas from NLP, e.g., Matpower \cite{5491276}, tend to 
produce some ``point'' for more instances than general-purpose NLP solvers, 
but the point can be far from feasible. 
Specialised solvers for the problem based on convexifications, which have been popular recently \cite{Low2014a,Low2014b}, 
 may produce excellent bounds on the objective function value, when the instance is feasible,
but both of the issues of detecting infeasibility of the instance and issues of 
the infeasibility of the point produced for a feasible instance remain.

In part, this is understandable, given the intrinsic complexity of the two related problems:
Testing feasibility is NP-Hard even for very simple cases 
 \cite{Lehmann2015},
  and assumed to be outside of NP, generally.
Testing feasibility of a strong, commonly-used convexification (SDP) relaxation is in NP $\cap$ co-NP \cite{Ramana1997},
similar to integer factorisation and (under mild assumptions) graph isomorphism,
which does not coincide with P, unless P = NP.
Also, arithmetic computation (AC) is reducible to the feasibility of SDPs \cite{TARASOV20082070}. 
Overall, this calls for a different approach. 

We propose an alternative, three-stage approach to polynomial optimisation, which compensates for the well-known shortcomings of local NLP solvers and convexifications in a principled fashion. In the three stages, we:
\begin{enumerate}
	\item find the closest feasible instance, by the projection onto the set of feasible inputs
	\item compute the optimum for the closest feasible instance to a limited precision
	\item find the closest feasible solution to the output of stage 2,
	by projection onto the set of feasible solutions.
\end{enumerate}

In terms of mathematical optimisation, the first stage considers an inner approximation of a 
polynomial optimisation problem,
possibly doubling the run-time of considering the ACOPF as a polynomial optimisation problem \cite{Las01} on its own. 
The third stage 
considers an outer approximation, but 
can be implemented using an adaption of a test of whether running Newton method for zero-finding  
 starting from an approximate solution is guaranteed to converge to the nearest zero,
 which is computationally rather efficient.

Our contributions are:
\begin{itemize}
\item the three-stage approach;
\item a semidefinite programming (SDP) relaxation for the projection onto the set of feasible inputs
aiming to develop a robust test of feasibility;
\item a numerical study of the projections onto the set of feasible inputs, comparing the SDP relaxation against NLP solvers;
\item a numerical study of projections onto the feasible set corresponding to the feasible input
 in order to obtain feasible solutions, which are based on our recent work \cite{liu2017hybrid}.
\end{itemize}


\section{A Three-Stage Approach}

Let us consider a polynomial optimisation problem (POP) parametrised by a vector $a$:
\begin{equation}\label{pp}
\begin{array}{ll}
\min_{x} & ~ f_0(x, a) \tag{POP} \\
\mathrm{s.t.}~& ~f_k(x, a)\geq 0,~k=1,...,m, \notag
\end{array}
\end{equation}
where the objective function $f_0$ and the constraints $f_k$, for $k=1,...,m$, 
are defined by multi-variate polynomials $f(x) = \sum_{\alpha \in \mathbb{N}^l : \lvert \alpha \rvert \leq \text{deg}_f} c_\alpha (f) x^\alpha$, with degree $\text{deg}_f$
in a decision variable $x \in \R^n$ 
and the parameter vector $a \in \R^v$.
Let $X(a)$ be the semi-algebraic feasible set for the decision variable $x$, 
as defined by parameters $a$.
Let us consider norms $\ell_p$, e.g., $p = \infty, 1, 2$.



We suggest a three-stage approach to POP:
\begin{algorithm}[h]
	\caption{A Three-Stage Approach to POP} 
	\label{algo:POP}
	\begin{algorithmic}[1]
		\STATE Project onto the set of feasible inputs:
		\begin{align*}
		\displaystyle\min_{b \in \R^v} & \displaystyle || a - b ||_p \label{s1POP} \text{ s.t. } X(b) \not = \emptyset 
		\end{align*}
		\STATE Assuming $X(b)$ has a proper interior, optimise over $X(b)$ to obtain a bound on the objective function value over $X(b)$, and the point $\tilde x$ at which this objective function is achieved, not necessarily within $X(b)$.
		\STATE Project $\tilde x$ onto the set $X(b)$:
		\begin{align*}
		\displaystyle\min_{x\in \R^n} & ||x - \tilde{x}||_p \text{ s.t. }\ x \in X(b)
		\end{align*}		
		\RETURN $x \in X(b)$ that is close to $\tilde{x}$. 
	\end{algorithmic}
\end{algorithm}

Within Line 2, one can utilise  semidefinite programming (SDP) relaxations of POP \cite{parrilo2000structured,Las01,parrilo2003semidefinite,anjos2011handbook}. Specifically, given a non-negative integer relaxation order $\omega \leq \omega_{\text{max}}$, where $\omega_{\text{max}} = \text{max} \{ \lceil \frac{1}{2} \text{deg}_{f_k} \rceil : k=0, \dots, m \}$, the dense SDP relaxation for POP of order $\omega$ due to \cite{Las01} is given by: 

\begin{equation}\label{sdp-pp}
\begin{array}{ll}
\min & ~ \sum_{\eta \in \text{supp}(f_0) } c_\eta (f_0, a) y_\eta \tag{SDP-POP} \\
\mathrm{s.t.}~& ~ M_{\omega - \omega_k} (f_k, a, y)\succeq 0,~k=1,...,m, \notag\\
& ~ M_{\omega} (y)\succeq 0, \notag
\end{array}
\end{equation}

where $\text{supp}(f_0)$ denotes the support of a function,  $M_\omega(y)$ and $M_\omega(f, y)$ denote the moment matrix and the localizing matrix for polynomial $f$ of order $\omega$, respectively. We refer to \cite{Las01,mevissen2010sdp,anjos2011handbook} for details. Notably, it is known that if POP satisfies conditions that are slightly stronger than the compactness of its feasible set, then the optimal value of (SDP-PP) converges to the one of POP for $\omega \longrightarrow \infty$.

For the projection onto the feasible inputs, we can construct inner approximations, 
e.g., as suggested by Henrion and Louembet \cite{Henrion2012}, who 
obtain the inner approximations by iteratively
 minimising curvature along algebraic varieties defining the boundary of \eqref{pp}. 
When the minimal curvature is negative, one obtains a separating hyperplane. 
There are also substantially more sophisticated approaches, e.g., \cite{6097030,6606873},
but these come at a substantial cost.

For the projection onto the feasible set, 
one could use Newton method, if one were close-enough to a local optimum.
A test, whether given a system
of polynomial
equalities and inequalities, and optimisation problems over such systems (POP), has been developed by Liu et al. \cite{liu2017hybrid}. 
This is based on an earlier work on a test, as to whether given a system of (complex-valued) polynomial equations and
a point, Newton method for zero-finding converges, which has been developed by Smale \cite{smale1986algorithms}, Shub \cite{ShubSmale1993}, 
Cucker \cite{Cucker1999}, and others \cite{DEREN1995,alphacertified} on the interface of algebraic geometry and numerical analysis, 
sometimes known as $\alpha$-$\beta$ theory.
Such a test whether one is close-enough to an optimum could become the termination criterion for the solver of Line 2.


More formally, the long history of work \cite{smale1986algorithms,ShubSmale1993,Cucker1999,DEREN1995} can be summarised as follows:
Consider a general real-valued polynomial system $f:\R^m \mapsto \R^n$, i.e.,
 a system of polynomial equations $f := ( f_1, \, \dots, \, f_n)$ in variables $x := (x_1, \dots, x_m) \in \R^m$.
Consider Newton operator at $x \in \R^m$ as 
$$N_f(x) := x - [\nabla f(x)]^\dagger f(x),$$
where $[\nabla f(x)]^\dagger \in \R^{m \times n}$ is the Moore-Penrose inverse of the Jacobian matrix of $f$ at $x$.
A sequence with initial point $x_0$ and iterates of Newton method subsequently, $x_{i + 1} := N_f(x_i)$ for $i \ge 0$,
is well-defined if $[\nabla f(x_i)]^\dagger$ is well defined at all $x_i, i \ge 0$.
We say that $x \in \R^m$ is an \emph{approximate zero} of $f$ if and only if:
\begin{enumerate}
    \item the sequence $\{x_i\}$ is well-defined
    \item there exists $x' \in \R^m$ such that $f(x') = 0$ and for all $i \geq 0$: $$\| x_i - x' \| \leq (1/2)^{2^{i-1} - 1} \|x_0 - x'\|.$$
\end{enumerate}
We call $x' \in \R^m$ the \emph{associated zero} of $x \in \R^m$ and 
say that $x$ \emph{represents} $x'$.
A key result of $\alpha$-$\beta$ theory is:

\begin{proposition}[Shub and Smale, \cite{ShubSmale1993}]\label{prop:alphabeta}
Let $f:\R^m \mapsto \R^n$ be a system of polynomial equations and define functions $\alpha(f,x), \beta(f,x), \gamma(f,x)$ as:

\begin{subequations}
\begin{align}
    \alpha(f,x) &:= \beta(f,x) \gamma(f,x), \label{alpha} \\
    \beta(f,x) &:= \left\| [\nabla f(x)]^\dagger  f(x)  \right\| = \| x - N_f(x) \|, \label{beta} \\
    \gamma(f,x) &:= \sup_{k > 1} \left\| {\frac{{ [\nabla f(x)]^\dagger [\nabla^{(k)}f] (x)}}{{k!}}} \right\|^{1/(k-1)}, \label{gamma}
\end{align}
\end{subequations}
where $[\nabla f(x)]^\dagger \in \R^{m \times n}$ is the Moore–-Penrose inverse of the Jacobian matrix of $f$ at $x$ and $[\nabla^{(k)}f] $ is the symmetric tensor whose entries are the $k$-th partial derivatives of $f$ at $x$.
Then there is a universal constant $\alpha_0 \in \R$ such that if $\alpha(f, x) \leq \alpha_0$, then $x$ is an approximate zero of $f$.
Moreover, if $x'$ denotes its associated zero, then $\| x - x'\| \le 2 \beta(f,x)$.
It can be shown that $\alpha_0 = \frac{13 - 3 \sqrt{17}}{4} \approx 0.157671$ satisfies this property.
\label{alphabeta}
\end{proposition}

The work of Liu et al. \cite{liu2017hybrid} can be summarised as:

\begin{proposition}[Liu et al. \cite{liu2017hybrid}]
There exists 
a universal constant $\alpha_0 \in \R$,
such that
for every instance of POP, 
there exists $\delta \in \R, \delta \ge 0$ 
and a function $\alpha: \R^m \mapsto \R$
specific to the instance of POP, such that for any 
$\epsilon > \delta$ and 
vector $x \in \R^m$ if $\alpha(x) \le \alpha_0$, then $x$ is 
in the domain of monotonicity of an optimum of the instance of POP, which is no more than 
$\epsilon$ away from the value of the global optimum with respect to its objective function.
\end{proposition}

While we refer to \cite{liu2017hybrid} for the complete details,
we stress that the test has been implemented and is practical.

\section{An Illustration}
\label{sec:prob}

Let us illustrate the three-stage approach on the example of ACOPF \cite{molzahn_hiskens-fnt2019}, the prototypical problem in power systems optimization. 
There, the complicated structure of the feasible set is understood, for certain small instances \cite[e.g.]{7879340}.

We consider the same representation of a power system as used by \cite{lavaei2012zero,Molzahn2011,GhaMM16,7879340} and 
the corresponding notation.  
The power system is represented by a directed graph, where each vertex $k \in N$ is called a ``bus'' and 
each directed edge $(l, m) \in E \subseteq N \times N$ is called a ``branch''. 
Each branch has an ideal phase-shifting transformer at its ``from'' end, and is modelled as a $\Pi$-equivalent circuit. 
Let $G \subseteq N$ be the set of slack buses, typically a singleton. 
Let $L \subseteq E$ be the remaining generators.
Let the remainder $N \setminus G$ represent the demands.

The corresponding constants are:

\begin{supertabular}{p{2cm}p{6cm}}
	$y \in \mathbb{R}^{|N|\times|N|}$ & network admittance matrix\\
	$\bar{b}_{lm}$             & shunt element value at branch $(l,m) \in E $\\
	$g_{lm}+jb_{lm}$           & series admittance on a branch $(l,m) \in E$\\
	$P^{d}_k$                  & active load (demand) at bus $k \in N$\\
	$Q^{d}_k$                  & reactive load (demand) at bus $k \in N$\\
	$P^{d}$                    &  aggregate active demand on period $t$\\
	$c_k^2, c_k^1, c_k^0, $    & coefficients of the quadratic generation costs $C_k$ at generator $k$\\
	$P_k^{\min}$, $P_k^{\max}$ & limits on active generation at bus $k$\\
	$Q_k^{\min}$, $Q_k^{\max}$ & limits on reactive generation at bus $k$\\
	$V_k^{\min}$, $V_k^{\max}$ & limits on the absolute value of the voltage at bus $k  \in N$\\
	$S_{lm}^{\max}$            & limit on the absolute value of the apparent power of branch $(l,m) \in L$,\\
	$e_k$                      & $k^{th}$ standard basis vector in $\mathbb{R}^{|N|}$ \\
\end{supertabular}

\vspace{0.2cm}

with the associated power-flow matrices:
\begin{align*}
y_k&=e_ke_k^T y, \\
y_{lm}&=(j\frac{\bar{b}_{lm}}{2}+g_{lm}+jb_{lm})e_le_l^T - (g_{lm}+jb_{lm})e_le_m^T, \\
Y_k&= \frac{1}{2} \left[ \begin{matrix} \Re(y_k  + y_k^T)  &\Im(y_k^T -y_k)  \\
\Im(y_k  - y_k^T))& \Re(y_k  + y_k^T) \end{matrix}  \right], \\
\bar{Y}_k&= -\frac{1}{2} \left[ \begin{matrix} \Im(y_k  + y_k^T)& \Re(y_k -y_k^T) \\
\Re(y_k^T  - y_k) & \Im(y_k  + y_k^T) \end{matrix}  \right], \\
M_k&=  \left[ \begin{matrix}e_ke_k^T & 0  \\
0 & e_ke_k^T \end{matrix}  \right], \\
Y_{lm}&= \frac{1}{2} \left[ \begin{matrix} \Re(y_{lm}  + y_{lm}^T)& \Im(y_{lm}^T -y_{lm}) \\
\Im(y_{lm} - y_{lm}^T) & \Re(y_{lm}  + y_{lm}^T) \end{matrix}  \right],\\
\bar{Y}_{lm}&=- \frac{1}{2} \left[ \begin{matrix} \Im(y_{lm}  + y_{lm}^T)& \Re(y_{lm}^T -y_{lm}) \\
\Re(y_{lm}^T - y_{lm}) & \Im(y_{lm}  + y_{lm}^T) \end{matrix}  \right].
\end{align*}


Using the usual rectangular power-voltage formulation of power flows in each period, the decision variables are:
\begin{supertabular}{lp{0.5\linewidth}}
	$x=\{\Re{V_{k}}+ j \Im{V_{k}}\}_{k \in N}$ & vector of voltages $V_{k}$,\\
	$ (P^{g}_k, Q^{g}_k)$ & active and reactive power of the generator at bus $k \in N$,\\
	$ (P_{lm}, Q_{lm})$ & active and reactive power flow on $(l,m) \in E$.\\
\end{supertabular}

We can hence formulate the polynomial optimisation problem (POP) of degree $2$, referred to as [$OP_2$] in \cite{GhaMM16}:

\begin{align}
\displaystyle\min & \displaystyle \sum_{k \in G} \bigg( c_k^2(P_k^g)^2+c_k^1 P_k^g+c_k^0 \bigg)& \label{objACOPF}\\
\text{s.t.}\ & P_k^{\min} \leq P_k^g \leq P_k^{\max} & \label{limact}\\
& Q_k^{\min} \leq Q_k^g \leq Q_k^{\max} & \label{limreact}\\
& P_k^g=\text{tr}(Y_kxx^T)+P_k^d& \label{defPkg}\\
& Q_k^g=\text{tr}(\bar{Y}_kxx^T)+Q_k^d& \label{defQkg}\\
& (V_k^{\min})^2 \leq \text{tr}(M_kxx^T) \leq (V_k^{\max})^2 & \label{limvolt}\\
& (P_{lm})^2+(Q_{lm})^2\leq (S_{lm}^{\max})^2  & \label{limthermOP2}\\
& P_{lm}=\text{tr}(Y_{lm}xx^T) & \label{defPlm}\\
& Q_{lm}=\text{tr}(\bar{Y}_{lm}xx^T) & \label{defQlm}
\end{align}
Objective function \eqref{objACOPF} is the cost of power generation. Constraints \eqref{limact} and \eqref{limreact} impose 
a bound on the active and reactive power, respectively. Constraints \eqref{defPkg} express the relationship between 
the active power, the admittance matrix, and the active load. Similarly, \eqref{defQkg} deal with reactive power. 
Constraints \eqref{limvolt} restrict the voltage on a given bus. 
Constraints \eqref{limthermOP2}, \eqref{defPlm} and \eqref{defQlm} limit the apparent power flow at each end of a given line.
In the following sections, we detail an approach for the projections onto the related sets.

A well-known convexification of the POP is given by the semidefinite programming (SDP) relaxation, which is referred to as Optimisation 3 by  \cite{lavaei2012zero}.
Optimisation 3 of \cite{lavaei2012zero} is obtained from $OP_2$ by replacing $xx^T$ in $OP_2$ with a matrix $W \succeq 0$, and relaxing the condition $\text{rank}(W) = 1$. 
By Theorem 1 of \cite{GhaMM16}, this coincides with the first relaxation order of \eqref{sdp-pp}, the moment-SOS hierarchy \cite{parrilo2000structured,Las01,parrilo2003semidefinite}.
In particular, the SDP relaxation reads:
\begin{align}
\displaystyle\min & \displaystyle \sum_{k \in G} \alpha_k & \label{objACOPF-SDP}\\
\text{s.t.}\ & P_k^{\min} \leq P_k^g \leq P_k^{\max} & \label{limact-SDP}\\
& Q_k^{\min} \leq Q_k^g \leq Q_k^{\max} & \label{limreact-SDP}\\
& P_k^g=\text{tr}(Y_k W)+P_k^d& \label{defPkg-SDP}\\
& Q_k^g=\text{tr}(\bar{Y}_k W)+Q_k^d& \label{defQkg-SDP}\\
& (V_k^{\min})^2 \leq \text{tr}(M_k W) \leq (V_k^{\max})^2 & \label{limvolt-SDP}\\
& \begin{bmatrix}
-(S_{lm}^{\max}) & \text{tr}(Y_{lm} W) & \text{tr}(\bar{Y}_{lm} W)  \\
\text{tr}(Y_{lm} W) & -1 & 0 \\
\text{tr}(\bar{Y}_{lm} W) & 0 & -1 \\
\end{bmatrix} \succeq 0  & \label{limthermOP2-SDP}\\
& \begin{bmatrix}
c_k^1 \text{tr}(Y_k W) - \alpha_k +a_k & \sqrt{c_k^2} \text{tr}(Y_k W) + b_k\\
\sqrt{c_k^2} \text{tr}(Y_k W) + b_k & -1
\end{bmatrix} \succeq 0 & \label{obj-Schur-SDP}\\
& W \succeq 0 & \label{W-SDP}\\
& \alpha_k \geq 0 & \label{alpha-SDP}
\end{align}
where $a_k = c_k^0 + c_k^1 P_k^d$, $b_k = \sqrt{c_k^2} P_k^d$, and $\alpha_k, k \in G$ are decision variables introduced to express the objective function via the Schur's complement formula.


\begin{algorithm}[tbh]
	\caption{A Three-Stage Approach for ACOPF} 
	\label{algo:two-stage}
	\begin{algorithmic}[1]
		\STATE Parametrised by slacks  $s$ of \eqref{acopfslacks}, let $X(s)$ be the feasible set of constraints \eqref{slack_act_up}-\eqref{slack_volt_down}, \eqref{defPkg}, \eqref{defQkg}, \eqref{limthermOP2}, \eqref{defPlm}, \eqref{defQlm}.
		\STATE Let $\chi = (x, \{P_k^g, Q_k^g\}_{k \in N}, \{P_{lm}, Q_{lm}\}_{(l,m) \in E})$.
		\STATE \label{Stage1} Stage 1: Set up projection $S_1$ onto the set of feasible inputs:
		\begin{align*}
		\displaystyle\min & \displaystyle ||s||_p \\
		\text{s.t.}\ & \chi \in X(s) & \\
		& s \geq 0 & 
		\end{align*}
		\STATE Run solver on $S_1$ to obtain upper bound $UB_1$ on the optimal value of $S_1$ and lower bound $LB_1$ on the optimal value of $S_1$.
		\STATE \label{Stage2} Stage 2: Set up feasible set $X(s)$, and the corresponding optimisation problem $S_2$:
		\begin{align*}
		\displaystyle\min & \ \eqref{objACOPF} \\
		\text{s.t.}\ & \chi \in X(s) \\
		& ||s||_p \leq UB_1 & \\
		& s \geq 0. & 
		\end{align*}	
		\STATE \label{Stage2} Run solver on $S_2$ to obtain point $\tilde{\chi}$, not necessarily within $X$, which is close to an optimum within $X$. 
		\STATE \label{Stage3} Stage 3: Set up a projection $S_3(\tilde{\chi})$ of $\tilde{\chi}$ onto the set $X(s)$ of feasible solutions:
		\begin{align*}
		\displaystyle\min_{\chi} & ||\chi - \tilde{\chi}||_p \\
		\text{s.t.}\ & \eqref{limact}-\eqref{defQlm}.
		\end{align*}		
		\RETURN Approximate solution of $S_3(\tilde{\chi})$, i.e., a feasible point $\chi$ that is close to $\tilde{\chi}$. 
	\end{algorithmic}
\end{algorithm}


In practical terms, within ACOPF, the active load $P^{d}_k$ and, to a lesser extent, the reactive load 
$Q^{d}_k$ are time-varying. 
With each update of $P^{d}_k, Q^{d}_k$, we woud like to decide whether it allows for a
feasible solution, or requires a redispatch, which may involve
increasing the limit $P_k^{\max}$ on active power generation, which would correspond to international transfers
or reserves, or
allowing for higher limits on absolute values $S_{lm}^{\max}$ of the apparent powers,
when the ambient conditions (temperature, wind-speed) prevent the temperature of the respective 
branch from increasing dangerously.

In mathematical terms, this corresponds to the question whether particular
 choices of parameters ($P^{d}_k$, $Q^{d}_k$, etc) make the 
POP (\ref{objACOPF}-\ref{defQlm}) feasible. 
Clearly, when one encounters an input, whose feasibility is not guaranteed,
one could project it onto the set of feasible inputs and report the differences.
To do so, we introduce slack variables $s \in \mathbb{R}^{\lvert G \rvert + \lvert G \rvert + \lvert N \rvert}$ obtained by concatenation of:
\begin{align}
 [ s_{P, k}^+, s_{P, k}^-, s_{Q, k}^+, s_{Q, k}^- ]_{k \in G} \textrm{ and } [s_{V, k}^+, s_{V, k}^-]_{k \in N}  \label{acopfslacks}\end{align}
for some fixed ordering of $G, N$,
wherein non-negative components quantify the extent to
 which a constraint is infeasible for formulation $OP_2$ (\ref{objACOPF}-\ref{defQlm}). 
Given an inequality constraint with $\leq$ sign, such as the upper bound expressed 
 in \eqref{limact}, a slack variable $s^+ \geq 0$  measures the amount by which the left-hand 
 side is greater than the right-hand side. 
Conversely, in inequalities with $\geq$ sign, slack $s^- \geq 0$ is the amount by which 
the left-hand side is lower than the right-hand side. 
 Hence, for the active power bounds \eqref{limact}, the slacks $\{ s_{P, k}^+, s_{P, k}^-,\}_{k \in G}$ are introduced as follows:
\begin{align}
 & P_k^g - s_{P, k}^+ \leq P_k^{\max} & \label{slack_act_up}\\
 & P_k^g + s_{P, k}^- \geq P_k^{\min} & \label{slack_act_down}
\end{align}
Slack variables $\{ s_{Q, k}^+, s_{Q, k}^-,\}_{k \in G}\cup \{ s_{V, k}^+, s_{V, k}^-,\}_{k \in N}$ can be introduced in an analogous manner for reactive power bounds \eqref{limreact} 
and voltage bounds \eqref{limvolt}: 
\begin{align}
& Q_k^g - s_{Q, k}^+ \leq Q_k^{\max} & \label{slack_react_up}\\
& Q_k^g + s_{Q, k}^- \geq Q_k^{\min} & \label{slack_react_down}\\
& \text{tr}(M_kxx^T) - s_{V, k}^+ \leq (V_k^{\max})^2 & \label{slack_volt_up} \\
& \text{tr}(M_kxx^T) + s_{V, k}^- \geq (V_k^{\min})^2 & \label{slack_volt_down} 
\end{align}

In practice, one can use even a limited-precision approximations of the lower and upper 
bounds ($LB_1$ and $UB_1$) on the constraint violations measured in terms of the slack variables. 
If $LB_1$ is found to be strictly non-negative, then the OPF instance can be declared as infeasible. 
Otherwise, in a second stage, one can amend the bounds on the slacks and search for an OPF solution 
$\tilde{\chi}$ minimizing \eqref{objACOPF}: this is described in Step \ref{Stage2} of Algorithm \ref{algo:two-stage}. In an  analogous manner, one can consider the SDP relaxation  (\ref{objACOPF-SDP}--
\ref{alpha-SDP}) to set up the three-stage approach for a relaxation of ACOPF, so as to cope with the infeasibility of active and reactive power bounds, and voltage limitations, i.e., constraints \eqref{limact-SDP}, \eqref{limreact-SDP}, and  \eqref{limvolt-SDP}.



To illustrate the computational performance, we have tested the approach
 on several IEEE test systems modified so as to exhibit infeasibilities. 
The modifications applied to the original IEEE test systems are summarised 
in Table \ref{table:instances}. 
Two implementations of our approach are described in the following subsections.

\subsection{Non-linear Programming (NLP)}

To illustrate the performance of non-linear programming solvers using the 
approach, we implemented the first two stages of the solution approach via AMPL models \cite{ampl}, 
which calls the non-linear programming solvers such as Ipopt \cite{wachter2006implementation}.
The communication between the stages is orchestrated in a Python framework. 
Finally, the Stage 3 is implemented using the Newton refinement step in PYPOWER \cite{Pypower}, 
which is a Python reimplementation of the MATPOWER package \cite{5491276}. 

\begin{table}[th]
	\centering
	\begin{tabular}{c|c}
		Name & Mod\\\hline
		case9-P70 & $P^{max}$ lowered by $70\%$,\\ & $P^{min}$ increased by $70\%$  \\\hline
		case14-P70 & $P^{max}$ lowered by $70\%$,\\ & $P^{min}$ increased by $70\%$\\\hline
		case14-V40 & [$V^{max}$, $V^{min}$] restricted by $40\%$ \\\hline
		case14-Q-80 & $Q^{max}$ lowered by $70\%$,\\ & $Q^{min}$ increased by $70\%$ \\ \hline
		case118-P60 & $P^{max}$ lowered by $0\%$,\\ & $P^{min}$ increased by $60\%$\\\hline
	\end{tabular}
	\caption{Characteristics of the instances tested with the two-stage approach.}
	\label{table:instances}
\end{table}


We have initialised Ipopt runs with the MATPOWER solution found on the unperturbed instances,
to aid the convergence within $10000$ iterations. 
We remark that Ipopt searches for local optima, and hence the solution obtained at Stage $1$ is
 a valid upper bound at Stage $2$. 
 
On case9-P70, Ipopt hits the maximum number of iterations at Stage 2. 
 Stage 3 corrects the Stage 2 solution by $8\%$. 
 On case14-Q-80, Ipopt converges to a locally infeasible solution and then Stage 3 
 incurs in numerical failure. On case118-P60, the PYPOWER refinement fails as well.
For all the instances tested apart from case118-P60, the slack variables are strictly non-negative only for 
 the constraints that are violated by the perturbations. This demonstrates that the slacks are correctly 
 identifying the infeasibility, if the $\ell_1$ norm is used to measure their impact.

A preliminary testing with the Global Optimisation (GO) solver Couenne \cite{belotti2009couenne} 
failed to converge to optimal solutions for both stages $1$ and $2$, within $120$ seconds of computation. 
The lower bounds $LB_1$ found for our instances were never strictly greater than $0$, 
therefore the GO solver did not provide certificates of infeasibility. 
Since the best solutions found by Couenne were comparable, 
or worse, to the Ipopt solutions, we only report the computational results obtained with Ipopt in Tables \ref{table:results_1} 
and \ref{table:results_inf}, for $\ell_1$ and $\ell_\infty$ norms, respectively.


\begin{table}[th]
	\centering
		\begin{tabular}{c|ccc|c}
		Instances	& \multicolumn{3}{c|}{Solution Values} & Time ($s$)\\
			& $S_1$ & $S_2$	& $S_3$	& 	\\\hline
		case9-P70 & $0.71$ & $5853.84$ & $5438.32$ & $127.93$ \\
		case14-P70 & $0.30$ & $8579.01$ & $8171.73$ & $15.37$  \\
		case14-V40 & $0.06$ & $10082.46$ & $8171.73$ &  $38.73$ \\
		case14-Q80 & $0.01$ & $14652.14*$ & - & $6.41$   \\
		case118-P60 & $55.39$ & $199196.31$ & $-$ & $527.61$  \\
		\hline
	\end{tabular}
	\caption{Solution values and computational times of Ipopt, obtained for the three stages. The $\ell_1$ norm is used to measure the usage of slack variables. The $*$ indicates a locally infeasible solution, and $-$ a failure on Stage 3.}
	\label{table:results_1}
\end{table}

\begin{table}[th]
	\centering
	\begin{tabular}{c|ccc|c}
		Instances	& \multicolumn{3}{c|}{Solution Values} & Time ($s$)\\
		& $S_1$ & $S_2$	& $S_3$	& 	\\\hline
		case9-P70 & $0.24$ & $5343.52$ & $5438.32$ & $9.38$ \\
		case14-P70 & $0.06$ & $8533.08^*$ & $8171.73$ & $120.47$  \\
		case14-V40 & $0.01$ & $14617.42^*$ & - &  $5.50$ \\
		case14-Q80 & $0.01$ & $10081.07$ & $8171.73$ & $55.73$ \\
		case118-P60 & $0.88$ & $147797.15$ & $-$ & $776.81$ \\
		\hline
	\end{tabular}
	\caption{Solution values and computational times of Ipopt, obtained for the three stages. The $\ell_\infty$ norm is used to measure the usage of slack variables. The $*$ indicates a locally infeasible solution, and $-$ a failure on Stage 3.}
	\label{table:results_inf}
\end{table}

%
When $\infty_\infty$ norm is used, if the PyPower Newton method converges, then it obtains the same solution of the $\ell_1$ norm. The bounds on the slacks are considerably smaller than those obtained for $\ell_1$ norm, and this contributes to finding Stage $2$ of better cost. In each instance, all slack variables are activated, regardless of the cause of the infeasibility.

\subsection{Semidefinite Programming (SDP)}

To illustrate the performance of our approach using convexifications, we have
 implemented stages $1$ and $2$ of our approach in the MATLAB toolbox YALMIP \cite{1393890}. 
In particular, we considered the semidefinite programming relaxation (\ref{objACOPF-SDP}--
\ref{alpha-SDP}) described in Section \ref{sec:prob}.
Although there are many SDP solvers \cite{marevcek2017low}, Table \eqref{table:resultsSDPnorm1} 
 reports the results obtained by using the SDP solver SeDuMi \cite{sturm1999using} and the $\ell_1$ norm 
 to measure the slack variables. The Newton method of Stage $3$ is then run in MATPOWER. 
 
\begin{table}[th]
	\centering
	\begin{tabular}{c|c|c|c|c}
		Instances	& \multicolumn{3}{c|}{Solution Values} & Time ($s$)\\
		& $S_1$ & $S_2$	& $S_3$ &  \\\hline
		case9-P70 & $0.70$ & $3759.97$ & $5438.32$ & $0.93$\\
		case14-P70 &  $0.30$ & $7853.13$ &  $-$  & $1.22$\\
		case14-V40 &  $0.01$ & $9881.35$ & $-$  & $1.32$ \\
		case14-Q80 & $0.06$ & $9789.38$ & $-$ & $1.30$ \\
		case118-P60 &  $/$ & $/$ & $/$ & $/$ \\
		\hline
	\end{tabular}
	\caption{Solution values and computational times of SeDuMi, obtained for the three stages. The $\ell_1$ norm is used to measure the usage of slack variables. The $/$ indicates that the solver ran out of memory, and $-$ a failure on Stage 3.}
	\label{table:resultsSDPnorm1}
\end{table}

\begin{table}[th]
	\centering
	\begin{tabular}{c|c|c|c|c}
		Instances	& \multicolumn{3}{c|}{Solution Values} & Time ($s$)\\
		& $S_1$ & $S_2$	& $S_3$ &  \\\hline
		case9-P70 & $0.079$ & $4247.01$ & $5438.32$ & $0.93$\\
		case14-P70 &  $0.022$ & $7753.23$ &  $-$  & $1.48$\\
		case14-V40 &  $0.003$ & $9719.31$ & $-$  & $1.69$ \\
		case14-Q80 & $0.004$ & $9767.04$ & $-$ & $1.60$ \\
		case118-P60 &  $/$ & $/$ & $/$ & $/$ \\
		\hline
	\end{tabular}
	\caption{Solution values and computational times of SeDuMi, obtained for the three stages. 
		The $\ell_\infty$ norm is used to measure the usage of slack variables. 
		The $/$ indicates that the solver ran out of memory, and $-$ a failure at Stage 3.}
	\label{table:resultsSDPnorminf}
\end{table}



%

The bounds on the slack variables, found at Stage 1, are not considerably smaller than 
those found using the POP formulation. 
The refined solution of case9-P70 obtained using the SDP relaxation coincides with the 
one obtained from the POP. 
However, case118-P60 is not solved by SeDuMi, because of memory limitations. 
As for the solutions found on the POP formulation, the slacks are activated for the constraints 
violated by the perturbations. On all instances obtained by perturbing the case14 instance, 
the Newton method fails to converge in $100$ iterations.
By solving an SDP relaxation of the POP, the three-stage approach is considerably faster. 
While the Stage 2 SDP solutions for case9-P70 and case14-Q80 are far from the Ipopt solutions 
by $35.77\%$ and $33.19\%$, the deviations for case14-P70 and case14-V40 are of $8.46\%$ and $1.99\%$, respectively. 
Table \ref{table:resultsSDPnorminf} displays the results obtained by using the $\ell_\infty$ norm on 
the slack variables. As observed for the POP solution, when the $\ell_\infty$ norm is used, 
the activation of the slacks is no longer localised in the constraints causing the infeasibility. 
With respect to the $\ell_1$ norm solutions, the Stage $1$ solution values are considerably smaller. 






\section{Conclusions and Future Work}

We proposed a three-stage approach for dealing with infeasibility 
both in the instances on the input and the output of a solver
for transmission-constrained 
problems in the alternating-current model. 
The approach compensates for well-known issues
in NLP solvers and convexifications at the price of solving an additional POP or its restriction.
Of independent interest could be the fact that slack variables make it possible to quantify the infeasibility. 
Our numerical evaluation shows that the introduction of slack variables in the POP formulation $OP_2$ penalised by 
the $\ell_1$ norm makes it possible to identify which constraints cause the infeasibilities. 
The approach is simple to implement and practical results can be obtained in short run-time. 
This could be of considerable interest to power systems practitioners.
\\

The method can plausibly be adapted to other families non-linear optimisation problems,
and perhaps analysed at that level of generality.
One may also wonder how to extend this to the on-line setting \cite{8442544}.
There is hence a considerable scope for further work.

\section*{Acknowledgment}
We thank Jie Liu, Alan Claude  Liddell, and Martin Tak{\'a}{\v{c}} for discussions related to these topics.

\addcontentsline{toc}{section}{Acknowledgment}


\bibliographystyle{IEEEtran}
\bibliography{two_stage,tcuc,enclosing}

\begin{thebibliography}{10}
\providecommand{\url}[1]{#1}
\csname url@samestyle\endcsname
\providecommand{\newblock}{\relax}
\providecommand{\bibinfo}[2]{#2}
\providecommand{\BIBentrySTDinterwordspacing}{\spaceskip=0pt\relax}
\providecommand{\BIBentryALTinterwordstretchfactor}{4}
\providecommand{\BIBentryALTinterwordspacing}{\spaceskip=\fontdimen2\font plus
\BIBentryALTinterwordstretchfactor\fontdimen3\font minus
  \fontdimen4\font\relax}
\providecommand{\BIBforeignlanguage}[2]{{%
\expandafter\ifx\csname l@#1\endcsname\relax
\typeout{** WARNING: IEEEtran.bst: No hyphenation pattern has been}%
\typeout{** loaded for the language `#1'. Using the pattern for}%
\typeout{** the default language instead.}%
\else
\language=\csname l@#1\endcsname
\fi
#2}}
\providecommand{\BIBdecl}{\relax}
\BIBdecl

\bibitem{castillo2013computational}
\BIBentryALTinterwordspacing
A.~Castillo and R.~P. O’Neill, ``Computational performance of solution
  techniques applied to the {ACOPF},'' 2013. [Online]. Available:
  \url{http://www.ferc.gov/industries/electric/indus-act/market-planning/opf-papers/acopf-5-computational-testing.pdf}
\BIBentrySTDinterwordspacing

\bibitem{carpentier1962contribution}
J.~Carpentier, ``Contribution a l’etude du dispatching economique,''
  \emph{Bulletin de la Societe Francaise des Electriciens}, vol.~3, no.~1, pp.
  431--447, 1962.

\bibitem{5491276}
R.~D. Zimmerman, C.~E. Murillo-Sanchez, and R.~J. Thomas, ``Matpower:
  Steady-state operations, planning, and analysis tools for power systems
  research and education,'' \emph{IEEE Transactions on Power Systems}, vol.~26,
  no.~1, pp. 12--19, Feb 2011.

\bibitem{Low2014a}
S.~Low, ``Convex relaxation of optimal power flow -- part i: Formulations and
  equivalence,'' \emph{Control of Network Systems, IEEE Transactions on},
  vol.~1, no.~1, pp. 15--27, March 2014.

\bibitem{Low2014b}
------, ``Convex relaxation of optimal power flow -- part ii: Exactness,''
  \emph{Control of Network Systems, IEEE Transactions on}, vol.~1, no.~2, pp.
  177--189, June 2014.

\bibitem{Lehmann2015}
K.~Lehmann, A.~Grastien, and P.~V. Hentenryck, ``{AC}-feasibility on tree
  networks is {NP-Hard},'' \emph{IEEE Transactions on Power Systems}, vol.~31,
  no.~1, pp. 798--801, Jan 2016.

\bibitem{Ramana1997}
M.~V. Ramana, ``An exact duality theory for semidefinite programming and its
  complexity implications,'' \emph{Mathematical Programming}, vol.~77, no.~1,
  pp. 129--162, Apr 1997.

\bibitem{TARASOV20082070}
\BIBentryALTinterwordspacing
S.~P. Tarasov and M.~N. Vyalyi, ``Semidefinite programming and arithmetic
  circuit evaluation,'' \emph{Discrete Applied Mathematics}, vol. 156, no.~11,
  pp. 2070 -- 2078, 2008, in Memory of Leonid Khachiyan (1952 - 2005 ).
  [Online]. Available:
  \url{http://www.sciencedirect.com/science/article/pii/S0166218X07001370}
\BIBentrySTDinterwordspacing

\bibitem{Las01}
J.~B. Lasserre, ``Global optimization with polynomials and the problem of
  moments,'' \emph{SIAM Journal on Optimization}, vol.~11, no.~3, pp. 796--817,
  2001.

\bibitem{liu2017hybrid}
J.~Liu, A.~C. Liddell, J.~Mare{\v{c}}ek, and M.~Tak{\'a}{\v{c}}, ``Hybrid
  methods in solving alternating-current optimal power flows,'' \emph{IEEE
  Transactions on Smart Grid}, vol.~8, no.~6, pp. 2988--2998, 2017.

\bibitem{parrilo2000structured}
P.~A. Parrilo, ``Structured semidefinite programs and semialgebraic geometry
  methods in robustness and optimization,'' Ph.D. dissertation, California
  Institute of Technology, 2000.

\bibitem{parrilo2003semidefinite}
------, ``Semidefinite programming relaxations for semialgebraic problems,''
  \emph{Mathematical programming}, vol.~96, no.~2, pp. 293--320, 2003.

\bibitem{anjos2011handbook}
M.~F. Anjos and J.~B. Lasserre, \emph{Handbook on semidefinite, conic and
  polynomial optimization}.\hskip 1em plus 0.5em minus 0.4em\relax Springer
  Science \& Business Media, 2011, vol. 166.

\bibitem{mevissen2010sdp}
M.~Mevissen and M.~Kojima, ``Sdp relaxations for quadratic optimization
  problems derived from polynomial optimization problems,'' \emph{Asia-Pacific
  Journal of Operational Research}, vol.~27, no.~01, pp. 15--38, 2010.

\bibitem{Henrion2012}
D.~Henrion and C.~Louembet, ``Convex inner approximations of nonconvex
  semialgebraic sets applied to fixed-order controller design,''
  \emph{International Journal of Control}, vol.~85, no.~8, pp. 1083--1092,
  2012.

\bibitem{6097030}
D.~Henrion and J.~Lasserre, ``Inner approximations for polynomial matrix
  inequalities and robust stability regions,'' \emph{IEEE Transactions on
  Automatic Control}, vol.~57, no.~6, pp. 1456--1467, June 2012.

\bibitem{6606873}
D.~Henrion and M.~Korda, ``Convex computation of the region of attraction of
  polynomial control systems,'' \emph{IEEE Transactions on Automatic Control},
  vol.~59, no.~2, pp. 297--312, Feb 2014.

\bibitem{smale1986algorithms}
S.~Smale, ``Algorithms for solving equations,'' in \emph{Proceedings of the
  International Congress of Mathematicians, Vol. 1, 2 (Berkeley, Calif.,
  1986)}.\hskip 1em plus 0.5em minus 0.4em\relax Providence, RI: Amer. Math.
  Soc., 1987, pp. 172--195.

\bibitem{ShubSmale1993}
M.~Shub and S.~Smale, ``On the complexity of {B}ezout's theorem i --
  {G}eometric aspects,'' \emph{Journal of the American Mathematical Society},
  vol.~6, no.~2, 1993.

\bibitem{Cucker1999}
\BIBentryALTinterwordspacing
F.~Cucker and S.~Smale, ``Complexity estimates depending on condition and
  round-off error,'' \emph{J. ACM}, vol.~46, no.~1, pp. 113--184, Jan. 1999.
  [Online]. Available: \url{http://doi.acm.org/10.1145/300515.300519}
\BIBentrySTDinterwordspacing

\bibitem{DEREN1995}
W.~Deren and Z.~Fengguang, ``The theory of {Smale}'s point estimation and its
  applications,'' \emph{Journal of Computational and Applied Mathematics},
  vol.~60, no.~1, pp. 253 -- 269, 1995.

\bibitem{alphacertified}
J.~D. Hauenstein and F.~Sottile, ``Alphacertified: Certifying solutions to
  polynomial systems,'' \emph{ACM Transactions on Mathematical Software},
  vol.~38, no.~4, pp. 28:1--28:20, 2012.

\bibitem{molzahn_hiskens-fnt2019}
D.~K. Molzahn and I.~A. Hiskens, ``{A Survey of Relaxations and Approximations
  of the Power Flow Equations},'' \emph{Foundations and Trends in Electric
  Energy Systems}, vol.~4, no. 1-2, pp. 1--221, February 2019.

\bibitem{7879340}
D.~K. {Molzahn}, ``Computing the feasible spaces of optimal power flow
  problems,'' \emph{IEEE Transactions on Power Systems}, vol.~32, no.~6, pp.
  4752--4763, Nov 2017.

\bibitem{lavaei2012zero}
J.~Lavaei and S.~H. Low, ``Zero duality gap in optimal power flow problem,''
  \emph{IEEE Transactions on Power Systems}, vol.~27, no.~1, pp. 92--107, Feb
  2012.

\bibitem{Molzahn2011}
D.~K. Molzahn, J.~T. Holzer, B.~C. Lesieutre, and C.~L. DeMarco,
  ``Implementation of a large-scale optimal power flow solver based on
  semidefinite programming,'' \emph{IEEE Transactions on Power Systems},
  vol.~28, no.~4, pp. 3987--3998, Nov 2013.

\bibitem{GhaMM16}
B.~Ghaddar, J.~Marecek, and M.~Mevissen, ``Optimal power flow as a polynomial
  optimization problem,'' \emph{IEEE Transactions on Power Systems}, vol.~31,
  no.~1, pp. 539--546, Jan 2016.

\bibitem{ampl}
R.~Fourer, D.~M. Gay, and B.~W. Kernighan, \emph{AMPL: A Mathematical
  Programing Language}, S.~W. Wallace, Ed.\hskip 1em plus 0.5em minus
  0.4em\relax Berlin, Heidelberg: Springer Berlin Heidelberg, 1989.

\bibitem{wachter2006implementation}
A.~W{\"a}chter and L.~T. Biegler, ``On the implementation of an interior-point
  filter line-search algorithm for large-scale nonlinear programming,''
  \emph{Mathematical programming}, vol. 106, no.~1, pp. 25--57, 2006.

\bibitem{Pypower}
R.~Lincoln, ``Pypower,'' \url{https://github.com/rwl/PYPOWER}.

\bibitem{belotti2009couenne}
P.~Belotti, ``Couenne: a user’s manual,'' Lehigh University, Tech. Rep.,
  2009.

\bibitem{1393890}
J.~Lofberg, ``Yalmip : a toolbox for modeling and optimization in matlab,'' in
  \emph{2004 IEEE International Conference on Robotics and Automation (IEEE
  Cat. No.04CH37508)}, Sept 2004, pp. 284--289.

\bibitem{marevcek2017low}
J.~Mare{\v{c}}ek and M.~Tak{\'a}{\v{c}}, ``A low-rank coordinate-descent
  algorithm for semidefinite programming relaxations of optimal power flow,''
  \emph{Optimization Methods and Software}, vol.~32, no.~4, pp. 849--871, 2017.

\bibitem{sturm1999using}
J.~F. Sturm, ``Using {SeDuMi} 1.02, a {MATLAB} toolbox for optimization over
  symmetric cones,'' \emph{Optimization methods and software}, vol.~11, no.
  1-4, pp. 625--653, 1999.

\bibitem{8442544}
J.~Liu, J.~Marecek, A.~Simonetta, and M.~Takac, ``A coordinate-descent
  algorithm for tracking solutions in time-varying optimal power flows,'' in
  \emph{2018 Power Systems Computation Conference (PSCC)}, June 2018, pp. 1--7.

\end{thebibliography}

\end{document}